# On the regularity of the generalised golden ratio function

Simon Baker[1] and Wolfgang Steiner[2]

[1]*Mathematics Institute,*
*University of Warwick,*
*Coventry, CV4 7AL, UK*

[2]*IRIF, CNRS UMR 8243,*
*Université Paris Diderot - Paris 7,*
*Case 7014, 75205 Paris Cedex 13, France*

*Email contacts*:
`simonbaker412@gmail.com`
`steiner@liafa.univ-paris-diderot.fr`

September 5, 2016

## Abstract

Given a finite set of real numbers $A$, the generalised golden ratio is the unique real number $\mathcal{G}(A) > 1$ for which we only have trivial unique expansions in smaller bases, and have non-trivial unique expansions in larger bases. We show that $\mathcal{G}(A)$ varies continuously with the alphabet $A$ (of fixed size). What is more, we demonstrate that as we vary a single parameter $m$ within $A$, the generalised golden ratio function may behave like $m^{1/h}$ for any positive integer $h$. These results follow from a detailed study of $\mathcal{G}(A)$ for ternary alphabets, building upon the work of Komornik, Lai, and Pedicini (2011). We give a new proof of their main result, that is we explicitly calculate the function $\mathcal{G}(\{0,1,m\})$. (For a ternary alphabet, it may be assumed without loss of generality that $A = \{0,1,m\}$ with $m \in (1,2)$].) We also study the set of $m \in (1,2]$ for which $\mathcal{G}(\{0,1,m\}) = 1 + \sqrt{m}$, we prove that this set is uncountable and has Hausdorff dimension 0. We show that the function mapping $m$ to $\mathcal{G}(\{0,1,m\})$ is of bounded variation yet has unbounded derivative. Finally, we show that it is possible to have unique expansions as well as points with precisely two expansions at the generalised golden ratio.

*Mathematics Subject Classification* 2010: 11A63, 28A80.

*Key words and phrases*: Beta expansions, Generalised golden ratios, Combinatorics on words

## 1 Introduction and statement of results

Let $A := \{a_0, a_1, \ldots, a_d\}$ be a set of real numbers satisfying $a_0 < a_1 < \cdots < a_d$. We call $A$ an alphabet. Given $\beta > 1$ and $x \in \mathbb{R}$, we say that a sequence $(u_k)_{k=1}^\infty \in A^\mathbb{N}$ is a $\beta$-expansion for $x$ over the alphabet $A$ if

$$x = \sum_{k=1}^\infty \frac{u_k}{\beta^k}.$$

When the underlying alphabet is obvious we may simply refer to $(u_k)$ as a $\beta$-expansion. Expansions in non-integer bases were introduced by Rényi [8]. Perhaps the most well studied case is when $\beta \in (1,2]$ and $A = \{0,1\}$. For $\beta \in (1,2]$ and this choice of alphabet, $x$ has a $\beta$-expansion over $A$ if and only if $x \in [0, \frac{1}{\beta-1}]$. Moreover, a result of Erdős, Joó, and Komornik [3] states that if $\beta \in (1, \frac{1+\sqrt{5}}{2})$ then every $x \in (0, \frac{1}{\beta-1})$ has a continuum of $\beta$-expansions. This result is complemented by a theorem of Daróczy and Katai [2] which states that if $\beta \in (\frac{1+\sqrt{5}}{2}, 2]$ then there exists $x \in (0, \frac{1}{\beta-1})$ with a unique $\beta$-expansion. Note that the end points of the interval $[0, \frac{1}{\beta-1}]$ trivially have a unique $\beta$-expansion for any $\beta \in (1,2]$. The above demonstrates that the golden ratio acts as a natural boundary between the possible cardinalities the set of expansions can take. It is natural to ask whether such a boundary exists for more general alphabets.

Before we state the definition of a generalised golden ratio it is necessary to define the univoque set. Given an alphabet $A$ and $\beta > 1$ we set

$$\mathcal{U}_\beta(A) := \left\{(u_k)_{k=1}^\infty \in A^\mathbb{N} : \sum_{k=1}^\infty \frac{u_k}{\beta^k} \text{ has a unique } \beta\text{-expansion}\right\}.$$

We call $\mathcal{U}_\beta(A)$ the univoque set. Note that for any $\beta > 1$ and alphabet $A = \{a_0, \ldots, a_d\}$ satisfying $a_0 < a_1 < \cdots < a_d$, the points

$$\sum_{k=1}^\infty \frac{a_0}{\beta^k} \text{ and } \sum_{k=1}^\infty \frac{a_d}{\beta^k}$$

both have a unique expansion, so $\overline{a_0}$ and $\overline{a_d}$ are always contained in the univoque set. Here and throughout $\overline{w}$ denotes the infinite periodic word with period $w$. We are now in a position to define a generalised golden ratio for an arbitrary alphabet. Given an alphabet $A$, we call $\mathcal{G}(A) \in (1, \infty)$ the generalised golden ratio for $A$ if whenever $\beta \in (1, \mathcal{G}(A))$ we have $\mathcal{U}_\beta(A) = \{\overline{a_0}, \overline{a_d}\}$, and if $\beta > \mathcal{G}(A)$ then $\mathcal{U}_\beta(A)$ contains a non-trivial element.

Komornik, Lai, and Pedicini [4] were the first authors to make a thorough study of generalised golden ratios over arbitrary alphabets. Importantly they proved that for any alphabet $A$ a generalised golden ratio exists. For ternary alphabets, they showed that the generalised golden ratio varies continously with the alphabet. We extend this result to alphabets of arbitrary size.

**Theorem 1.** *Let $\Delta_d := \{(a_0, a_1, \ldots, a_d) \in \mathbb{R}^{d+1} : a_0 < a_1 < \cdots < a_d\}, d \geq 1$. The map $(a_0, a_1, \ldots, a_d) \mapsto \mathcal{G}(\{a_0, a_1, \ldots, a_d\})$ is continuous on $\Delta_d$.*

We prove this theorem in Section 2. In the rest of the paper, we restrict our attention to ternary alphabets. Every ternary alphabet can be assumed to be of the form $A = \{0, 1, m\}$ for some $m > 1$ because shifting the alphabet and multiplying by a constant does not affect the generalised golden ratio. We thus set

$$\mathcal{G}(m) := \mathcal{G}(\{0, 1, m\}) \text{ and } \mathcal{U}_\beta(m) := \mathcal{U}_\beta(\{0, 1, m\}).$$

Moreover,

$$\mathcal{G}(m) = \mathcal{G}(\{0, 1, m\}) = \mathcal{G}(\{-m, 1-m, 0\}) = \mathcal{G}\left(\left\{\frac{m}{m-1}, 1, 0\right\}\right) = \mathcal{G}\left(\frac{m}{m-1}\right).$$

By the above $\mathcal{G}(m) = \mathcal{G}(\frac{m}{m-1})$ and we may therefore assume $m \in (1, 2]$. The authors of [4] considered $m \geq 2$, and their results read as follows in our setting.

**Theorem KLP.** *The function $\mathcal{G} : (1, 2] \to \mathbb{R}$ is continuous and satisfies*

$$2 \leq \mathcal{G}(m) \leq 1 + \sqrt{m}$$

*for all $m \in (1, 2]$. Moreover, the following statements hold.*

- $\mathcal{G}(m) = 2$ *for $m \in (1, 2]$ if and only if $m = \frac{2^k}{2^k - 1}$ for some positive integer $k$.*

- *The set $\mathfrak{M} := \{m \in (1, 2] : \mathcal{G}(m) = 1 + \sqrt{m}\}$ is a Cantor set, its largest element is $x^2 \approx 1.7548$ where $x \approx 1.3247$ is the smallest Pisot number.*

- *Each connected component $(m_1, m_2)$ of $(1, x^2) \setminus \mathfrak{M}$ has a point $\mu$ such that $\mathcal{G}$ is strictly decreasing on $[m_1, \mu]$ and strictly increasing on $[\mu, m_2]$; $\mathcal{G}$ is strictly increasing on $[x^2, 2]$.*

In Section 3, we reprove all these results, making some of the statements more explicit and simplifying several proofs. An approximation of the graph of $\mathcal{G}$ can be found in Figure 3.1. The function $\mathcal{G}$ is given by implicit equations on subintervals of $(1, 2]$, and it has the following unusual regularity properties.

**Theorem 2.** *The function $\mathcal{G} : (1, 2] \to \mathbb{R}$ is differentiable except on the set $\mathfrak{M}$ and on the countable set of points $\mu$ defined in Theorem KLP. Its derivative is unbounded, but its total variation is less than 2.*

We have the following result on the size of $\mathfrak{M}$.

**Theorem 3.** *The set $\mathfrak{M}$ is an uncountable Cantor set with Hausdorff dimension $0$.*

On certain intervals, the function $\mathcal{G}$ has the following simple form.

**Theorem 4.** *Let $h$ be a positive integer and $2^h \leq m \leq \left(1 + \sqrt{\frac{m}{m-1}}\right)^h$. Then we have*

$$\mathcal{G}(m) = \mathcal{G}\left(\frac{m}{m-1}\right) = m^{1/h}.$$

Note that if $m = 2^h$ then $m < \left(1 + \sqrt{\frac{m}{m-1}}\right)^h$, thus the set of $m$ defined within Theorem 4 is non-empty.

Moreover, we have the following result on the size of the set of expansions at the generalised golden ratio.

**Theorem 5.** *There exists $m \in (1,2]$ such that:*

- $\mathcal{U}_{\mathcal{G}(m)}(m)$ *contains non-trivial elements.*
- *When $\beta = \mathcal{G}(m)$ there exists an $x$ with precisely two $\beta$-expansions.*

Finally, we remark that the problem of calculating $\mathcal{G}(A)$ remains wide open. Only $\mathcal{G}(\{0, 1, \ldots, m\})$ has been calculated for any positive integer $m$ in [1].

## 2 Continuity of $\mathcal{G}(A)$

Before proving Theorem 1, we recall results of Pedicini [7, Proposition 2.1 and Theorem 3.1] on (unique) expansions in non-integer bases over arbitrary alphabets; see also [4, Theorem 2.2].

**Theorem P.** *Let $\beta \in (1, q(A)]$, with*

$$q(A) := 1 + \frac{a_d - a_0}{\max\{a_1 - a_0, a_2 - a_1, \ldots, a_d - a_{d-1}\}}.$$

*Every $x \in \left[\frac{a_0}{\beta - 1}, \frac{a_d}{\beta - 1}\right]$ has a $\beta$-expansion over $A$. We have $u_1 u_2 \cdots \in \mathcal{U}_\beta(A)$ if and only if, for all $i \geq 1$,*

$$\sum_{k=0}^{\infty} \frac{u_{i+k}}{\beta^k} < a_{j+1} + \frac{a_0}{\beta - 1} \qquad \text{when } u_i = a_j \neq a_d, \tag{2.1}$$

*and*

$$\sum_{k=0}^{\infty} \frac{u_{i+k}}{\beta^k} > a_{j-1} + \frac{a_d}{\beta - 1} \qquad \text{when } u_i = a_j \neq a_0. \tag{2.2}$$

*Remark* 2.1. The conditions (2.1) and (2.2) can be restated in terms of uniqueness regions $E_a$: Let

$$E_{a_0} = \left[\frac{a_0 \beta}{\beta - 1}, a_1 + \frac{a_0}{\beta - 1}\right), \ E_{a_j} = \left(a_{j-1} + \frac{a_d}{\beta - 1}, a_{j+1} + \frac{a_0}{\beta - 1}\right), \ 1 \leq j < d, \ E_{a_d} = \left(a_{d-1} + \frac{a_d}{\beta - 1}, \frac{a_d \beta}{\beta - 1}\right].$$

Then (2.1) and (2.2) hold if and only if $\sum_{k=0}^{\infty} \frac{u_{i+k}}{\beta^k} \in E_{u_i}$.

By the following lemma, it is sufficient to consider $\beta \leq q(A)$.

**Lemma 2.2.** *We have $\mathcal{G}(A) \leq q(A)$.*

*Proof.* If $\beta > 1 + \frac{a_d - a_0}{a_{j+1} - a_j}$ for some $0 \leq j < d$, then $a_j + \frac{a_d}{\beta - 1} < a_{j+1} + \frac{a_0}{\beta - 1}$ and thus $a_j \overline{a_d} \in \mathcal{U}_\beta(A)$. $\square$

*Remark* 2.3. This upper bound is attained for certain alphabets. For example, let $A = \{0, 1, 4, 5\}$. For $\beta = q(A) = 8/3$, the uniqueness regions are $E_0 = [0, 1)$, $E_1 = (3, 4)$, $E_4 = (4, 5)$ and $E_5 = (7, 8]$. If $\sum_{k=0}^{\infty} \frac{u_{i+k}}{\beta^k} \in E_{u_i}$, then $\sum_{k=0}^{\infty} \frac{u_{i+1+k}}{\beta^k} \in (E_{u_i} - u_i)\beta$; the latter intervals are $[0, 8/3)$, $(16/3, 8)$, $(0, 8/3)$ and $(16/3, 8]$ respectively. Therefore, the only unique expansions are $\overline{0}$ and $\overline{5}$.

*Proof of Theorem 1.* As $\mathcal{G}(A) = \mathcal{G}\left(\frac{A-a_0}{a_d-a_0}\right)$, we have $\mathcal{G}(\{a_0, a_1, \ldots, a_d\}) = \mathcal{G}(\iota \circ r(a_0, a_1, \ldots, a_d))$, with

$$r: \Delta_d \to \Delta'_d, \quad (a_0, a_1, \ldots, a_d) \mapsto \left(\frac{a_1 - a_0}{a_d - a_0}, \frac{a_2 - a_0}{a_d - a_0}, \ldots, \frac{a_{d-1} - a_0}{a_d - a_0}\right),$$

$$\iota: \Delta'_d \to \mathcal{P}(\mathbb{R}), \quad (a_1, a_2, \ldots, a_{d-1}) \mapsto \{0, a_1, a_2, \ldots, a_{d-1}, 1\},$$

and $\Delta'_d = \{(a_1, a_2 \ldots, a_{d-1}) \in \mathbb{R}^{d-1} : 0 < a_1 < a_2 < \cdots < a_{d-1} < 1\}$. As $r$ is continuous on $\Delta_d$, it is sufficient to prove that $\mathcal{G} \circ \iota$ is continuous on $\Delta'_d$.

Let $\mathbf{a} = (a_1, a_2, \ldots, a_{d-1}) \in \Delta'_d$ and $\varepsilon > 0$ arbitrary but fixed. We will show that $|\mathcal{G}(\iota(\mathbf{b})) - \mathcal{G}(\iota(\mathbf{a}))| \leq 3\varepsilon$ for all $\mathbf{b}$ in a neighbourhood of $\mathbf{a}$. Let first $X \subset \Delta'_d$ be a closed neighbourhood of $\mathbf{a}$ such that $|q(\iota(\mathbf{b})) - q(\iota(\mathbf{a}))| \leq \varepsilon$ for all $\mathbf{b} \in X$. (Note that $q \circ \iota$ is continuous on $\Delta'_d$.) Set

$$\alpha = \min_{\mathbf{b} \in X} q(\iota(\mathbf{b})) - \varepsilon, \qquad Y = \{\mathbf{b} \in X : \mathcal{G}(\iota(\mathbf{b})) < \alpha\}.$$

If $Y = \emptyset$, then $X$ is a neighbourhood of $\mathbf{a}$ with $|\mathcal{G}(\iota(\mathbf{b})) - \mathcal{G}(\iota(\mathbf{a}))| \leq 2\varepsilon$ for all $\mathbf{b} \in X$. Otherwise, let $\ell \geq 2$ be such that $\sum_{k=1}^{\ell} \alpha^{-k} \geq (\alpha + \varepsilon - 1)^{-1}$. Then

$$b_{j+1} - b_j \leq \frac{1}{q(\iota(\mathbf{b})) - 1} \leq \frac{1}{\alpha + \varepsilon - 1} \leq \sum_{k=1}^{\ell} \frac{1}{\alpha^k} \tag{2.3}$$

for all $(b_1, \ldots, b_{d-1}) \in Y$, $0 \leq j < d$, with $b_0 = 0$, $b_d = 1$. Set

$$\delta(\mathbf{a}, \mathbf{b}) = \min_{0 \leq j < d} \left((a_{j+1} - a_j) - (b_{j+1} - b_j)\right)$$

(with $b_0 = a_0 = 0$, $b_d = a_d = 1$), and let $Z \subset X$ be a neighbourhood of $\mathbf{a}$ such that

$$\frac{a_j}{(\alpha + \varepsilon)^k} - \frac{b_j}{\alpha^k} \leq \delta(\mathbf{a}, \mathbf{b}), \quad \frac{b_j}{(\alpha + \varepsilon)^k} - \frac{a_j}{\alpha^k} \leq \delta(\mathbf{a}, \mathbf{b}) \quad \text{for all } 1 \leq j \leq d,\ 1 \leq k \leq \ell, \tag{2.4}$$

$$\frac{1 - a_j}{(\alpha + \varepsilon)^k} - \frac{1 - b_j}{\alpha^k} \leq \delta(\mathbf{a}, \mathbf{b}), \quad \frac{1 - b_j}{(\alpha + \varepsilon)^k} - \frac{1 - a_j}{\alpha^k} \leq \delta(\mathbf{a}, \mathbf{b}) \quad \text{for all } 0 \leq j < d,\ 1 \leq k \leq \ell, \tag{2.5}$$

for all $\mathbf{b} = (b_1, \ldots, b_{d-1}) \in Z$. Note that $\delta(\mathbf{a}, \mathbf{b}) \leq 0$, thus we also have

$$\frac{a_j}{\alpha + \varepsilon} \leq \frac{b_j}{\alpha}, \quad \frac{1 - a_j}{\alpha + \varepsilon} \leq \frac{1 - b_j}{\alpha}, \quad \frac{b_j}{\alpha + \varepsilon} \leq \frac{a_j}{\alpha}, \quad \frac{1 - b_j}{\alpha + \varepsilon} \leq \frac{1 - a_j}{\alpha} \quad \text{for all } 0 \leq j \leq d. \tag{2.6}$$

For $\mathbf{b} \in Y \cap Z$ and $\beta \in (\mathcal{G}(\iota(\mathbf{b})), \alpha]$, choose $\mathbf{u} = u_1 u_2 \cdots \in \mathcal{U}_\beta(\iota(\mathbf{b}))$. Assume, w.l.o.g., that $u_1 u_2 \notin \{00, 11\}$. We show first that $\mathbf{u}$ does not contain $\ell$ consecutive zeros or ones. Indeed, suppose that $u_{i+1} = u_{i+2} = \cdots = u_{i+\ell} = 1$ for some $i \geq 1$; then we have

$$\sum_{k=1}^{\infty} \frac{u_{i+k}}{\beta^k} \geq \sum_{k=1}^{\ell} \frac{1}{\beta^k} \geq \sum_{k=1}^{\ell} \frac{1}{\alpha^k} \geq b_{j+1} - b_j$$

for all $0 \leq j < d$, hence $u_i = 1$ because of (2.1); recursively we would obtain that $u_{i-1} = \cdots = u_1 = 1$, contradicting that $u_1 u_2 \neq 11$. Similarly, $u_{i+1} = u_{i+2} = \cdots = u_{i+\ell} = 0$ implies that $\sum_{k=1}^{\infty} (1 - u_{i+k}) \beta^{-k} \geq b_j - b_{j-1}$ for all $1 \leq j \leq d$, hence $u_i = 0$ because of (2.2), eventually contradicting that $u_1 u_2 \neq 00$.

We define the sequence $\tilde{\mathbf{u}}$ via the relation $\tilde{u}_{i+k} = a_j$ if $u_{i+k} = b_j$. Let now $i \geq 1$. We have $\tilde{u}_{i+k}(\beta + \varepsilon)^{-k} \leq u_{i+k} \beta^{-k}$ for all $k \geq 1$ because (2.6) implies that $a_j(\beta + \varepsilon)^{-k} \leq b_j \beta^{-k}$ for all $0 \leq j \leq d$. Moreover, (2.4) and (2.6) give that

$$\frac{a_j}{(\beta + \varepsilon)^k} - \frac{b_j}{\beta^k} \leq \frac{a_j}{(\alpha + \varepsilon)^k} - \frac{b_j}{\alpha^k} \leq \delta(\mathbf{a}, \mathbf{b})$$

for all $1 \leq k \leq \ell$, $1 \leq j \leq d$. Since $u_{i+k} \neq 0$ for some $1 \leq k \leq \ell$, we have $\tilde{u}_{i+k}(\beta + \varepsilon)^{-k} \leq u_{i+k} \beta^{-k} + \delta(\mathbf{a}, \mathbf{b})$ for some $k \geq 1$. Using (2.1), we get

$$\sum_{k=1}^{\infty} \frac{\tilde{u}_{i+k}}{(\beta + \varepsilon)^k} \leq \sum_{k=1}^{\infty} \frac{u_{i+k}}{\beta^k} + \delta(\mathbf{a}, \mathbf{b}) < b_{j+1} - b_j + \delta(\mathbf{a}, \mathbf{b}) \leq a_{j+1} - a_j \quad \text{when } u_i = b_j \neq 1.$$

Similarly, we obtain from (2.2), (2.5) and (2.6) that

$$\sum_{k=1}^{\infty} \frac{1-\tilde{u}_{i+k}}{(\beta+\varepsilon)^k} \le \sum_{k=1}^{\infty} \frac{1-u_{i+k}}{\beta^k} + \delta(\mathbf{a},\mathbf{b}) < b_j - b_{j-1} + \delta(\mathbf{a},\mathbf{b}) \le a_j - a_{j-1} \quad \text{when } u_i = b_j \ne 0.$$

Therefore, we have $\tilde{\mathbf{u}} \in \mathcal{U}_{\beta+\varepsilon}(\iota(\mathbf{a}))$, thus $\mathcal{G}(\iota(\mathbf{a})) \le \mathcal{G}(\iota(\mathbf{b})) + \varepsilon$ for all $\mathbf{b} \in Y \cap Z$.

For $\mathbf{b} \in X \setminus Y$, recall that $\mathcal{G}(\iota(\mathbf{a})) \le q(\iota(\mathbf{a})) \le \alpha + 2\varepsilon \le \mathcal{G}(\iota(\mathbf{b})) + 2\varepsilon$. Similarly, we obtain for all $\mathbf{b} \in Z$ that $\mathcal{G}(\iota(\mathbf{b})) \le \mathcal{G}(\iota(\mathbf{a})) + \varepsilon$ when $\mathbf{a} \in Y$, $\mathcal{G}(\iota(\mathbf{b})) \le \mathcal{G}(\iota(\mathbf{a})) + 3\varepsilon$ when $\mathbf{a} \notin Y$. This gives that $|\mathcal{G}(\iota(\mathbf{b})) - \mathcal{G}(\iota(\mathbf{a}))| \le 3\varepsilon$ for all $\mathbf{b} \in Z$, thus $\mathcal{G} \circ \iota$ is continuous at $\mathbf{a}$. $\square$

## 3 Generalised golden ratios over ternary alphabets

### 3.1 Statements

Komornik, Lai and Pedicini [4] described the function $m \mapsto \mathcal{G}(m)$ on the interval $(1, 2]$. We provide more details for this function, in particular for the set

$$\mathfrak{M} := \{m \in (1, 2] : \mathcal{G}(m) = 1 + \sqrt{m}\}.$$

For $h \ge 0$, let $\tau_h$ be the substitution on the alphabet $\{0, 1\}$ defined by

$$\tau_h(0) = 0^{h+1}1, \qquad \tau_h(1) = 0^h 1,$$

and set $S = \{\tau_h : h \ge 0\}$. A (right) infinite word $\mathbf{u}$ is a *limit word* of a sequence of substitutions $(\sigma_n)_{n \ge 0}$ if there exist words $\mathbf{u}^{(n)}$ with $\mathbf{u}^{(0)} = \mathbf{u}$ and $\mathbf{u}^{(n)} = \sigma_n(\mathbf{u}^{(n+1)})$ for all $n \ge 0$. A sequence $(\sigma_n)_{n \ge 0} \in S^{\mathbb{N}}$ is *primitive* if $\sigma_n \ne \tau_0$ for infinitely many $n \ge 0$. A limit word of a primitive sequence in $S^{\mathbb{N}}$ starts with $\sigma_0 \sigma_1 \cdots \sigma_n(0)$ for all $n \ge 0$ and is therefore unique. If $\sigma_n = \tau_0$ for all $n \ge 0$, then $1^k 0\bar{1}$, $k \ge 0$, and $\bar{1}$ are limit words of $(\sigma_n)_{n \ge 0}$; we are only interested in $0\bar{1}$ and $\bar{1}$. Therefore, we define the following sets of limit words (or $S$-adic words), where $S^* = \bigcup_{n \ge 0} S^n$ denotes the set of finite products of substitutions in $S$:

$$\mathcal{S} = \mathcal{S}_\infty \cup \mathcal{S}_{0\bar{1}} \cup \mathcal{S}_{\bar{1}} \quad \text{with} \quad \mathcal{S}_{0\bar{1}} = \{\sigma(0\bar{1}) : \sigma \in S^*\}, \quad \mathcal{S}_{\bar{1}} = \{\sigma(\bar{1}) : \sigma \in S^*\},$$
$$\mathcal{S}_\infty = \{\mathbf{u} : \mathbf{u} \text{ is the limit word of a primitive sequence of substitutions in } S^{\mathbb{N}}\}.$$

*Remark* 3.1. Komornik, Lai and Pedicini [4] observed that the sequences $\mathbf{u} \in \mathcal{S}_\infty$ with the leading 0 removed are exactly the standard Sturmian sequences. However, they omitted the word "standard".

For $\mathbf{u} = u_0 u_1 \cdots \in \{0, 1\}^{\mathbb{N}}$, we define $\mathfrak{m}_\mathbf{u} \ge 1$ as the unique solution to

$$\mathfrak{m}_\mathbf{u} = 1 + \sum_{k=0}^{\infty} \frac{u_k}{(1 + \sqrt{\mathfrak{m}_\mathbf{u}})^k}. \tag{3.1}$$

*Remark* 3.2. We can rewrite (3.1) as

$$1 + \sqrt{\mathfrak{m}_\mathbf{u}} = 2 + \sum_{k=0}^{\infty} \frac{u_k}{(1 + \sqrt{\mathfrak{m}_\mathbf{u}})^{k+1}},$$

i.e., Parry's [6] $\beta$-expansion of $\beta = 1 + \sqrt{\mathfrak{m}_\mathbf{u}}$ is $2\mathbf{u}$. We have $\mathfrak{m}_\mathbf{u} = 1$ if and only if $\mathbf{u} = \bar{0}$.

For $\sigma \in S^*$, we define the interval $I_\sigma = [\mathfrak{m}_{\sigma(0\bar{1})}, \mathfrak{m}_{\sigma(\bar{1})}] \subset \left(1, \frac{3+\sqrt{5}}{2}\right]$. We define $\beta_\sigma \ge 2$ implicitly via the equation

$$1 + \sum_{k=1}^{\infty} \frac{\tilde{u}_k^{(\sigma)}}{\beta_\sigma^k} = (\beta_\sigma - 1)\left(1 + \sum_{k=0}^{\infty} \frac{u_k^{(\sigma)}}{\beta_\sigma^{k+1}}\right),$$

where

$$\tilde{u}_0^{(\sigma)} \tilde{u}_1^{(\sigma)} \tilde{u}_2^{(\sigma)} \cdots = \sigma(0\bar{1}), \qquad u_0^{(\sigma)} u_1^{(\sigma)} u_2^{(\sigma)} \cdots = \sigma(\bar{1}).$$

Moreover, we let $\mu_\sigma$ denote the coinciding value, i.e.,

$$\mu_\sigma := 1 + \sum_{k=1}^{\infty} \frac{\tilde{u}_k^{(\sigma)}}{\beta_\sigma^k}.$$

Note that all the numbers and sequences do not change if we replace $\sigma$ by $\sigma\tau_0$ since $\tau_0(0\bar{1}) = 0\bar{1}$ and $\tau_0(\bar{1}) = \bar{1}$. Therefore, we can assume that $\sigma \in S^* \setminus S^*\tau_0$.

The following propositions recover the main part of Theorem KLP, adding explicit equations giving the generalised golden ratios, which are used for drawing the graph of $\mathcal{G}$ in Figure 3.1.

**Proposition 3.3.** *The interval* $\left(1, \frac{3+\sqrt{5}}{2}\right]$ *admits the partition*
$$\{I_\sigma : \sigma \in S^* \setminus S^*\tau_0\} \cup \{\{\mathfrak{m}_\mathbf{u}\} : \mathbf{u} \in \mathcal{S}_\infty\}. \tag{3.2}$$

**Proposition 3.4.** *For $\sigma \in S^*$ and $m \in \left[\mathfrak{m}_{\sigma(0\bar{1})}, \mu_\sigma\right]$, $\mathcal{G}(m)$ is given by*
$$m = 1 + \sum_{k=1}^{\infty} \frac{\tilde{u}_k^{(\sigma)}}{\mathcal{G}(m)^k}.$$

*For $\sigma \in S^*$ with $\sigma(1) \neq 1$ and $m \in \left[\mu_\sigma, \mathfrak{m}_{\sigma(\bar{1})}\right]$, $\mathcal{G}(m)$ is given by*
$$\frac{m}{\mathcal{G}(m) - 1} = 1 + \sum_{k=1}^{\infty} \frac{u_k^{(\sigma)}}{\mathcal{G}(m)^{k+1}}.$$

**Proposition 3.5.** *We have*
$$\mathfrak{M} = \{\mathfrak{m}_\mathbf{u} : \mathbf{u} \in \mathcal{S} \setminus \{\bar{1}\}\}.$$

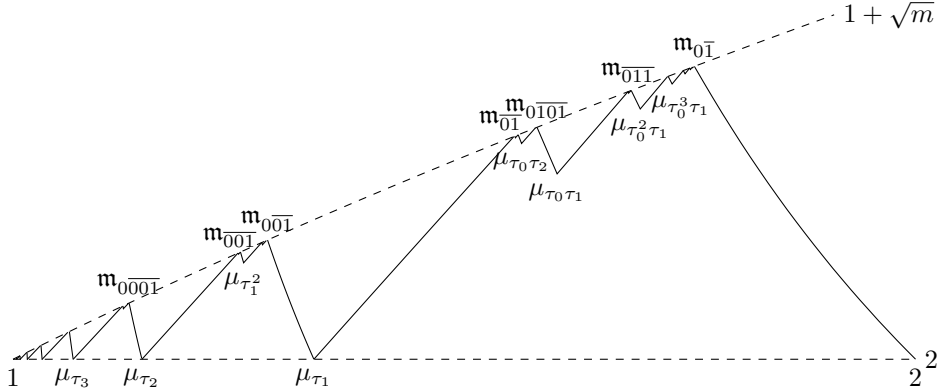

Figure 3.1: A graph of $\mathcal{G}(m)$.

## 3.2 Partition of $\left(1, \frac{3+\sqrt{5}}{2}\right]$

We first prove Proposition 3.3, using the following lemmas.

**Lemma 3.6.** *Let $\sigma \in S^*$. Then $\sigma$ preserves the lexicographic order on infinite words.*

*Proof.* The lexicographic order on infinite words is preserved by the identity and by $\sigma \in S$. By induction on $n$, this holds for all $\sigma \in S^n$, $n \geq 0$. $\square$

**Lemma 3.7.** *Let $\sigma \in S^*$ with $\sigma(1) \neq 1$, and write $\sigma(1) = 0w1$. Then*
$$\sigma(0\bar{1}) = 0\,\overline{w01}.$$
*In particular, $1w0$ is a circular shift of $\sigma(1)$.*

*Proof.* Since $\sigma(1) = \sigma\tau_0(1)$ and $\sigma(0\bar{1}) = \sigma\tau_0(0\bar{1})$, we assume w.l.o.g. that $\sigma = \tau_{h_0}\tau_{h_1}\cdots\tau_{h_n}$ with $n \geq 0$, $h_n \neq 0$. Let $\sigma_k = \tau_{h_k}$. Then $\sigma_{[0,n]}(1) = 0w1$ with
$$w = \sigma_0(1)\,\sigma_{[0,1]}(1)\cdots\sigma_{[0,n-1]}(1)\sigma_{[0,n-1]}(0^{h_n-1})\,\sigma_{[0,n-2]}(0^{h_{n-1}})\cdots\sigma_0(0^{h_1})0^{h_0}.$$

Let $v = \sigma_0(1)\,\sigma_{[0,1]}(1)\cdots\sigma_{[0,n-1]}(1)$. Then we have $\sigma_{[0,n]}(0) = 0w01v$ and $v\sigma_{[0,n]}(1) = w01v$. Therefore, $1w0$ is a circular shift of $\sigma(1)$ and $\sigma(0\bar{1}) = 0\,\overline{w01}$. $\square$

**Lemma 3.8.** *Let* $\mathbf{u} = u_0 u_1 \cdots \in \{0,1\}^{\mathbb{N}} \setminus \{\overline{0}\}$. *We have* $\mathbf{u} \in \mathcal{S}$ *if and only if*

$$u_0 u_1 u_2 \cdots \leq u_i u_{i+1} u_{i+2} \cdots \leq 1 u_1 u_2 \cdots \qquad \text{for all } i \geq 0. \tag{3.3}$$

*Proof.* Assume that (3.3) holds. Then $u_0 = 0$ or $\mathbf{u} = \overline{1} = \tau_0(\overline{1})$. If $u_0 = 0$, let $h \geq 0$ be minimal such that $u_{h+1} = 1$. Then each 1 is followed by $0^{h+1}1$ or $0^h1$, i.e., $\mathbf{u} = \tau_h(\mathbf{u}')$ for some word $\mathbf{u}' = u_0' u_1' \cdots$. Moreover, we have $\mathbf{u}' \leq u_i' u_{i+1}' \cdots \leq 1 u_1' u_2' \cdots$ for all $i \geq 0$. In case $\mathbf{u}' = \overline{0}$, we have $\mathbf{u} = \tau_{h+1}(\overline{1})$. Therefore, we can repeat the arguments and obtain recursively that $\mathbf{u}$ is the limit word of a sequence $(\sigma_n)_{n \geq 0} \in S^{\mathbb{N}}$. More precisely, we have $\mathbf{u} \in \mathcal{S}_{0\overline{1}}$ or $\mathbf{u}$ starts with $\sigma_{[0,n]}(0)$ for all $n \geq 0$, i.e., $\mathbf{u} \in \mathcal{S}_{\infty} \cup \mathcal{S}_{\overline{1}}$.

Consider now $\mathbf{u} \in \mathcal{S}_{\infty} \cup \mathcal{S}_{0\overline{1}}$, limit word of $(\sigma_n)_{n \geq 0} \in S^{\mathbb{N}}$. Then $\mathbf{u}$ starts with $\sigma_{[0,n]}(0)$ for all $n \geq 0$. Denote the preimage of $\mathbf{u}$ by $\sigma_0$ by $\mathbf{u}' = u_0' u_1' \cdots$, i.e., $\sigma_0(\mathbf{u}') = \mathbf{u}$. Suppose that $u_i u_{i+1} \cdots \leq \mathbf{u}$. Then $u_i u_{i+1} \cdots$ starts with $\sigma_0(0)$, and $u_i u_{i+1} \cdots = \sigma_0(u_{i'}' u_{i'+1}' \cdots)$ for some $i' \geq 0$. This implies that $u_{i'}' u_{i'+1}' \cdots \leq \mathbf{u}'$, thus $u_{i'}' u_{i'+1}' \cdots$ starts with $\sigma_1(0)$. Inductively, we obtain that $u_i u_{i+1} \cdots$ starts with $\sigma_{[0,n]}(0)$ for all $n \geq 0$, i.e., $u_i u_{i+1} \cdots = \mathbf{u}$. Suppose now that $u_i u_{i+1} \cdots \geq 1 u_1 u_2 \cdots$. Then $u_i = 1$ and $u_{i+1} u_{i+2} \cdots = \sigma_0(u_{i'}' u_{i'+1}' \cdots)$ for some $i' \geq 0$, with $u_{i'}' u_{i'+1}' \cdots \geq 1 u_1' u_2' \cdots$. We get that

$$u_i u_{i+1} \cdots = 1 \sigma_0(1) \sigma_{[0,1]}(1) \sigma_{[0,2]}(1) \cdots = 1 u_1 u_2 \cdots.$$

Therefore, (3.3) holds.

Finally, let $\mathbf{u} \in \mathcal{S}_{\overline{1}}$. If $\mathbf{u} = \overline{1}$, then (3.3) holds trivially. Otherwise, we have $\mathbf{u} = \sigma \tau_h(\overline{1})$ with $\sigma \in S^*$, $h \geq 1$. Then $\sigma \tau_{h-1} \tau_j(0\overline{1}) \in \mathcal{S}_{0\overline{1}}$ converges to $\mathbf{u}$ for $j \to \infty$ (in the usual topology of infinite words). By the previous paragraph, (3.3) holds for these words. Hence, it also holds for the limit word $\mathbf{u}$. $\square$

*Proof of Proposition 3.3.* For $m \in \left(1, \frac{3+\sqrt{5}}{2}\right]$, Parry's $(1+\sqrt{m})$-expansions of $1+\sqrt{m}$ are of the form $2\mathbf{u} \in \{0,1,2\}^{\mathbb{N}}$ with $\overline{0} < \mathbf{u} \leq \overline{1}$ and are ordered with respect to the lexicographic order. We show that the interval $(\overline{0}, \overline{1}] \subset \{0,1,2\}^{\mathbb{N}}$ admits the partition

$$\{[\sigma(0\overline{1}), \sigma(\overline{1})] : \sigma \in S^* \setminus S^* \tau_0\} \cup \{\{\mathbf{u}\} : \mathbf{u} \in \mathcal{S}_\infty\}.$$

Assume that $\mathbf{u} = u_0 u_1 \cdots \notin \mathcal{S}$, i.e., (3.3) does not hold. Let $i \geq 1$ be minimal such that one of the equalities is not satisfied. Suppose first that $u_i u_{i+1} \cdots < \mathbf{u}$; then $u_{i-1} = 1$. Similarly to the proof of Lemma 3.8, let $\sigma_0 \in S$ be such that $u_0 \cdots u_{i-1} = \sigma_0(u_0' \cdots u_{i'-1}')$ with $u_0' \cdots u_{i'-1}' \neq 0 \cdots 0$. By minimality of $i$, we have $u_{i'-1}' = 1$, and $u_0' \cdots u_{i'-1}' = 1 \cdots 1$ implies $i' = 1$. Therefore, we can define recursively substitutions $\sigma_j \in S$ until

$$\sigma_{[0,n]}(1) = u_0 u_1 \cdots u_{i-1}.$$

Then we have $\mathbf{u} < \sigma_{[0,n]}(1) \mathbf{u} < \cdots < \sigma_{[0,n]}(\overline{1})$. By Lemma 3.6, we have $\sigma_{[0,n]}(0\overline{1}) < \mathbf{u}$.

Suppose now that $u_i u_{i+1} \cdots > 1 u_1 u_2 \cdots$. Then we have substitutions $\sigma_k = \tau_{h_k}$ such that

$$\sigma_{[0,n]}(0) = u_0 u_1 \cdots u_{i-1} 1 \sigma_0(1) \sigma_{[0,1]}(1) \cdots \sigma_{[0,n-1]}(1),$$

with $h_n \neq 0$. We have $\mathbf{u} > u_0 \overline{u_1 \cdots u_{i-1} 1}$, and the latter word is equal to $\sigma_{[0,n]}(0\overline{1})$ by Lemma 3.7 and its proof. Since $u_0 \cdots u_{i-1} < \sigma_{[0,n]}(1)$, we also have $\mathbf{u} < \sigma_{[0,n]}(\overline{1})$.

We have seen that each $\mathbf{u}$ is the limit word of a primitive sequence of substitutions $\boldsymbol{\sigma} \in S^{\mathbb{N}}$ or between the extremal limit words of a non-primitive sequence $\boldsymbol{\sigma} \in S^{\mathbb{N}}$. To see that $\boldsymbol{\sigma}$ is unique, let $\mathbf{u}$ and $\tilde{\mathbf{u}}$ be limit words of two different sequences $(\sigma_n)_{n \geq 0}$ and $(\tilde{\sigma}_n)_{n \geq 0}$. Let $n \geq 0$ be minimal such that $\sigma_n \neq \tilde{\sigma}_n$. Let $\sigma_n = \tau_h$, $\tilde{\sigma}_n = \tau_j$, and assume w.l.o.g. that $h < j$. Then we have $\tilde{\mathbf{u}} \leq \tilde{\sigma}_{[0,n]}(\overline{1}) \leq \sigma_{[0,n]}(\overline{0}) < \mathbf{u}$. Therefore, the intervals are disjoint. $\square$

### 3.3 Calculating the generalised golden ratio

We now prove that $\mathcal{G}(m)$ is as in Theorem KLP and Proposition 3.4.

**Lemma 3.9.** *Let* $m \in (1, 2]$, $\beta \in [m, m+1]$, *and* $\mathbf{u} = u_0 u_1 \cdots \in \{0,1\}^{\mathbb{N}} \setminus \{\overline{0}\}$. *Then* $\mathbf{u} \in \mathcal{U}_\beta(m)$ *if and only if*

$$\frac{m}{\beta - 1} < 1 + \sum_{k=1}^{\infty} \frac{u_{i+k}}{\beta^k} < m \quad \text{for all } i \geq 0 \text{ such that } u_i = 1.$$

*Proof.* As $q(\{0,1,m\}) = 1 + m$, Theorem P and Remark 2.1 give that $\mathbf{u} \in \mathcal{U}_\beta(m)$ if and only if $\sum_{k=0}^\infty \frac{u_{i+k}}{\beta^k} \in E_{u_i}$ for all $i \geq 0$, with $E_0 = [0,1)$ and $E_1 = (\frac{m}{\beta-1}, m)$. If $u_i = 0$, then we either have $u_j = 0$ for all $j \geq i$ and thus $\sum_{k=0}^\infty \frac{u_{i+k}}{\beta^k} = 0 \in E_0$ or there exists $j > i$ such that $u_i = \cdots = u_{j-1} = 0$, $u_j = 1$. In the latter case, $\sum_{k=0}^\infty \frac{u_{j+k}}{\beta^k} \in E_1$ implies that $\sum_{k=0}^\infty \frac{u_{i+k}}{\beta^k} \in \beta^{i-j} E_1 \subset E_0$ since $\beta \geq m$. This proves the lemma. □

**Lemma 3.10.** *Let $\mathbf{u} \in \mathcal{S}_\infty$, $m \in (1,2]$. Then we have $\mathbf{u} \in \mathcal{U}_{1+\sqrt{m}}(m)$ if and only if $m = \mathfrak{m}_\mathbf{u}$. In particular, we have $\mathcal{G}(\mathfrak{m}_\mathbf{u}) \leq 1 + \sqrt{\mathfrak{m}_\mathbf{u}}$.*

*Proof.* By Lemma 3.9, we have $\mathbf{u} = u_0 u_1 \cdots \in \mathcal{U}_{1+\sqrt{m}}(m)$ if and only if

$$\sqrt{m} < 1 + \sum_{k=1}^\infty \frac{u_{i+k}}{(1+\sqrt{m})^k} < m$$

for all $i \geq 0$ such that $u_i = 1$. By Lemma 3.8 and since $\mathbf{u}$ is aperiodic, $u_i = 1$ implies that $\mathbf{u} < u_{i+1} u_{i+2} \cdots < u_1 u_2 \cdots$. Here, the bounds $\mathbf{u}$ and $u_1 u_2 \cdots$ cannot be improved because, for all $n \geq 0$, $1\sigma_{[0,n]}(0)$ and $1\sigma_0(1) \cdots \sigma_{[0,n-1]}(1)$ (which is a suffix of $\sigma_{[0,n]}(0)$) are factors of $\mathbf{u}$. Therefore, we have $\mathbf{u} \in \mathcal{U}_{1+\sqrt{m}}(m)$ if and only if

$$\sqrt{m} \leq 1 + \sum_{k=1}^\infty \frac{u_k}{(1+\sqrt{m})^{k+1}} \quad \text{and} \quad 1 + \sum_{k=1}^\infty \frac{u_k}{(1+\sqrt{m})^k} \leq m.$$

This means that $1 + \sum_{k=1}^\infty u_k (1+\sqrt{m})^{-k} = m$, i.e., $m = \mathfrak{m}_\mathbf{u}$. □

**Lemma 3.11.** *Let $\sigma \in S^*$ and $m > 1$. There is a unique number $f_\sigma(m) > 1$ such that*

$$m = 1 + \sum_{k=1}^\infty \frac{\tilde{u}_k^{(\sigma)}}{f_\sigma(m)^k}. \tag{3.4}$$

*We have $f_\sigma'(m) < 0$, $f_\sigma(\mathfrak{m}_{\sigma(0\bar{1})}) = 1 + \sqrt{\mathfrak{m}_{\sigma(0\bar{1})}}$, $f_\sigma(m) < 1 + \sqrt{m}$ if and only if $m > \mathfrak{m}_{\sigma(0\bar{1})}$, and $\sigma(\bar{1}) \notin \mathcal{U}_{f_\sigma(m)}(m)$ if $m \leq 2$.*

*Proof.* Let $h_m(x) = 1 + \sum_{k=1}^\infty \tilde{u}_k^{(\sigma)} x^{-k} - m$. Then $\lim_{x \to 1} h_m(x) = \infty$, $\lim_{x \to \infty} h_m(x) = 1 - m < 0$, $h_m(x)$ is continuous and strictly monotonically decreasing, thus $f_\sigma(m)$ is the unique solution of $h_m(x) = 0$. We have

$$\frac{1}{f_\sigma'(m)} = -\sum_{k=1}^\infty \frac{k \tilde{u}_k^{(\sigma)}}{f_\sigma(m)^{k+1}} < 0,$$

in particular $f_\sigma(m) < f_\sigma(\mathfrak{m}_{\sigma(0\bar{1})}) = 1 + \sqrt{\mathfrak{m}_{\sigma(0\bar{1})}} < 1 + \sqrt{m}$ for $m > \mathfrak{m}_{\sigma(0\bar{1})}$.

By Lemma 3.7, $1\tilde{u}_1^{(\sigma)} \tilde{u}_2^{(\sigma)} \cdots$ is a periodic word with the same period as $\sigma(\bar{1})$. Therefore, (3.4) and Lemma 3.9 imply that $\sigma(\bar{1}) \notin \mathcal{U}_{f_\sigma(m)}(m)$. □

**Lemma 3.12.** *Let $\sigma \in S^*$ and $m > 1$. There is a unique number $g_\sigma(m) > 1$ such that*

$$\frac{m}{g_\sigma(m) - 1} = 1 + \sum_{k=0}^\infty \frac{u_k^{(\sigma)}}{g_\sigma(m)^{k+1}}. \tag{3.5}$$

*We have $g_\sigma'(m) > 0$, $g_\sigma(\mathfrak{m}_{\sigma(\bar{1})}) = 1 + \sqrt{\mathfrak{m}_{\sigma(\bar{1})}}$, $g_\sigma(m) < 1 + \sqrt{m}$ if and only if $m < \mathfrak{m}_{\sigma(\bar{1})}$, and $\sigma(\bar{1}) \notin \mathcal{U}_{g_\sigma(m)}(m)$ if $m \leq 2$.*

*Proof.* Setting $h_m(x) = \frac{m}{x-1} - 1 - \sum_{k=0}^\infty \frac{u_k^{(\sigma)}}{x^{k+1}}$, we have

$$h_m'(x) = \sum_{k=0}^\infty \frac{(k+1) u_k^{(\sigma)}}{x^{k+2}} - \frac{m}{(x-1)^2} \leq \sum_{k=0}^\infty \frac{k+1}{x^{k+2}} - \frac{m}{(x-1)^2} = \frac{1-m}{(x-1)^2} < 0 \tag{3.6}$$

for $x > 1$. Similarly to the proof of Lemma 3.11, $g_\sigma(m)$ is the unique solution of $h_m(x) = 0$. (Note that $g_\sigma(m) = m$ if $\sigma(\bar{1}) = \bar{1}$.) We have

$$\frac{1}{g'_\sigma(m)} = \bigl(1 - g_\sigma(m)\bigr) h'_m(g_\sigma(m)) > 0.$$

Now, $g_\sigma(m) < 1 + \sqrt{m}$ is equivalent to $h_m(1+\sqrt{m}) < 0$, i.e., $1 + \sum_{k=0}^\infty u_k^{(\sigma)}(1+\sqrt{m})^{-k-1} > \sqrt{m}$. By Remark 3.2, we obtain that $m < \mathfrak{m}_{\sigma(\bar{1})}$. Since $1\sigma(\bar{1})$ is a suffix of $\sigma(\bar{1})$, (3.5) and Lemma 3.9 give that $\sigma(\bar{1}) \notin \mathcal{U}_{g_\sigma(m)}(m)$. $\square$

**Lemma 3.13.** *Let $\sigma \in S^*$. There is a unique $\mu_\sigma \in \bigl(\mathfrak{m}_{\sigma(0\bar{1})}, \mathfrak{m}_{\sigma(\bar{1})}\bigr)$ with $f_\sigma(\mu_\sigma) = g_\sigma(\mu_\sigma)$. We have $f_\sigma(\mu_\sigma) = g_\sigma(\mu_\sigma) \geq 2$, with equality if and only if $\sigma(\bar{1}) = \overline{0^n 1}$ for some $n \geq 0$.*

*If $m \in \bigl[\mathfrak{m}_{\sigma(0\bar{1})}, \mu_\sigma\bigr]$, then $\sigma(\bar{1}) \in \mathcal{U}_\beta(m)$ for all $\beta > f_\sigma(m)$, in particular $\mathcal{G}(m) \leq f_\sigma(m)$.*

*If $m \in \bigl[\mu_\sigma, \mathfrak{m}_{\sigma(\bar{1})}\bigr]$, then $\sigma(\bar{1}) \in \mathcal{U}_\beta(m)$ for all $\beta > g_\sigma(m)$, in particular $\mathcal{G}(m) \leq g_\sigma(m)$.*

*Proof.* The number $\mu_\sigma$ is well defined since $f'(m) < 0$, $g'(m) > 0$ on $I_\sigma$,

$$f_\sigma(\mathfrak{m}_{\sigma(0\bar{1})}) = 1 + \sqrt{\mathfrak{m}_{\sigma(0\bar{1})}} > g_\sigma(\mathfrak{m}_{\sigma(0\bar{1})}) \text{ and } f_\sigma(\mathfrak{m}_{\sigma(\bar{1})}) < 1 + \sqrt{\mathfrak{m}_{\sigma(\bar{1})}} = g_\sigma(\mathfrak{m}_{\sigma(\bar{1})}).$$

If $\sigma(\bar{1}) = \bar{1}$, then $\mu_\sigma = \beta_\sigma = 2$. Assume in the following that $\sigma(\bar{1}) \neq \bar{1}$ and let $m = 1 + \sum_{k=0}^\infty u_k^{(\sigma)} 2^{-k-1}$, i.e., $g_\sigma(m) = 2$. By Lemma 3.7, we have $\sigma(\bar{1}) = \overline{0w1} \leq \overline{w01} = \tilde{u}_1^{(\sigma)} \tilde{u}_2^{(\sigma)} \cdots$ for some finite word $w$, thus

$$1 + \sum_{k=1}^\infty \frac{\tilde{u}_k^{(\sigma)}}{2^k} \geq 1 + \sum_{k=0}^\infty \frac{u_k^{(\sigma)}}{2^{k+1}} = m,$$

hence $f_\sigma(m) \geq 2 = g_\sigma(m)$. This implies that $\beta_\sigma \geq 2$. If $\beta_\sigma = 2$, then we must have $0w = w0$, i.e., $w = 0 \cdots 0$. Therefore, $\beta_\sigma = 2$ is equivalent to $\sigma(\bar{1}) = \overline{0^n 1}$ for some $n \geq 0$.

Let now $m \in \bigl[\mathfrak{m}_{\sigma(0\bar{1})}, \mu_\sigma\bigr]$ and $\beta > f_\sigma(m)$, or $m \in \bigl[\mu_\sigma, \mathfrak{m}_{\sigma(\bar{1})}\bigr]$ and $\beta > g_\sigma(m)$. Then we also have $\beta > g_\sigma(m)$ and $\beta > f_\sigma(m)$ respectively. For $i \geq 0$ with $u_i^{(\sigma)} = 1$, we get

$$\frac{m}{\beta - 1} < 1 + \sum_{k=0}^\infty \frac{u_k^{(\sigma)}}{\beta^{k+1}} \leq 1 + \sum_{k=1}^\infty \frac{u_{i+k}^{(\sigma)}}{\beta^k} \leq 1 + \sum_{k=1}^\infty \frac{\tilde{u}_k^{(\sigma)}}{\beta^k} < m,$$

where the first inequality follows from $\beta > g_\sigma(m)$ and (3.6), the last inequality from $\beta > f_\sigma(m)$, and the middle inequalities are direct consequences of $\beta \geq 2$ and $\sigma(\bar{1}) \leq u_{i+1}^{(\sigma)} u_{i+2}^{(\sigma)} \cdots \leq \tilde{u}_1^{(\sigma)} \tilde{u}_2^{(\sigma)} \cdots$, which holds by Lemmas 3.8 and 3.7. Thus $\sigma(\bar{1}) \in \mathcal{U}_\beta(m)$. $\square$

The preceding lemmas show that $\mathcal{G}(m) \leq 1 + \sqrt{m}$ for all $m \in (1, 2]$. The next lemma justifies why we have restricted our attention to sequences in $\{0, 1\}^{\mathbb{N}}$.

**Lemma 3.14.** *Let $m \in (1, 2]$, $\beta \leq 1 + \sqrt{m}$, $u_0 u_1 \cdots \in \mathcal{U}_\beta(m)$. Then $u_i = m$ implies $u_0 \cdots u_i = m \cdots m$.*

*Proof.* If $u_i = m$, $i \geq 1$, then (2.2) implies that $u_{i-1} + \sum_{k=1}^\infty \frac{u_{i-1+k}}{\beta^k} > u_{i-1} + \frac{1}{\beta}\bigl(1 + \frac{m}{\beta-1}\bigr) \geq u_{i-1} + 1$, thus (2.1) excludes that $u_{i-1} = 0$ or $u_{i-1} = 1$. Recursively, we obtain that $u_k = m$ for all $0 \leq k \leq m$. $\square$

**Lemma 3.15.** *Let $m \in (1, 2]$, $\beta < 2$. Then $\mathcal{U}_\beta(m)$ is trivial.*

*Proof.* Let $u_0 u_1 \cdots \in \mathcal{U}_\beta(m)$. By Theorem P, we have $u_i \neq 1$ for all $i \geq 0$. Since $m \leq 1 + \frac{m}{\beta-1}$, we have $m\bar{0} \notin \mathcal{U}_\beta(m)$, thus Lemma 3.14 implies that $\mathcal{U}_\beta(m) = \{\bar{0}, \bar{m}\}$. $\square$

**Lemma 3.16.** *Let $m \in (1, 2]$, $\beta \leq 1 + \sqrt{m}$, and $u_0 u_1 \cdots \in \mathcal{U}_\beta(m) \cap (\{0,1\}^{\mathbb{N}} \setminus \{\bar{0}\})$. Then we have $\inf\{u_{i+1} u_{i+2} \cdots : i \geq 0, u_i = 1\} \in \mathcal{S}_\infty \cup \mathcal{S}_{\bar{1}}$.*

*Proof.* Let $\tilde{\mathbf{u}} = \tilde{u}_0 \tilde{u}_1 \cdots = \inf\{u_{i+1} u_{i+2} \cdots : i \geq 0, u_i = 1\}$. Since $\tilde{\mathbf{u}} = \bar{1} \in \mathcal{S}_{\bar{1}}$ when $\tilde{u}_0 = 1$, we assume in the following that $\tilde{u}_0 = 0$. For all $i \geq 0$ with $\tilde{u}_i = 1$, we have

$$\sum_{k=1}^\infty \frac{\tilde{u}_{i+k}}{\beta^k} < m - 1 \leq \beta\left(\frac{m}{\beta-1} - 1\right) \leq \beta \inf_{i \geq 0:\, u_i = 1} \sum_{k=1}^\infty \frac{u_{i+k}}{\beta^k} \leq \beta \sum_{k=0}^\infty \frac{\tilde{u}_k}{\beta^{k+1}} = \sum_{k=1}^\infty \frac{\tilde{u}_k}{\beta^k},$$

since $E_1 = \bigl(\frac{m}{\beta-1}, m\bigr)$ and $\beta \leq 1 + \sqrt{m}$. As $\beta \geq 2$ by Lemma 3.15, we obtain that $\tilde{u}_i \tilde{u}_{i+1} \cdots < 1 \tilde{u}_1 \tilde{u}_2 \cdots$ for all $i \geq 0$. By the definition of $\tilde{\mathbf{u}}$, we also have $\tilde{u}_i \tilde{u}_{i+1} \cdots \geq \tilde{\mathbf{u}}$, thus $\tilde{\mathbf{u}} \in \mathcal{S}$ by Lemma 3.8. Moreover, we have $\tilde{\mathbf{u}} \notin \mathcal{S}_{0\bar{1}}$ by Lemma 3.7 and the fact that $\tilde{u}_i \tilde{u}_{i+1} \cdots < 1 \tilde{u}_1 \tilde{u}_2 \cdots$ for all $i \geq 0$. $\square$

*Remark* 3.17. One obtains similarly that $\sup\{0u_{i+1}u_{i+2}\cdots : i \geq 0, u_i = 1\} \in \mathcal{S}_\infty \cup \mathcal{S}_{0\overline{1}}$.

**Proposition 3.18.** *We have*

$$\mathcal{G}(m) = \begin{cases} 1 + \sqrt{m} & \text{if } m \in \{\mathfrak{m}_{\mathbf{u}} : \mathbf{u} \in \mathcal{S}_\infty\}, \\ f_\sigma(m) & \text{if } m \in \left[\mathfrak{m}_{\sigma(0\overline{1})}, \mu_\sigma\right], \sigma \in S^*, \\ g_\sigma(m) & \text{if } m \in \left[\mu_\sigma, \mathfrak{m}_{\sigma(\overline{1})}\right], \sigma \in S^*. \end{cases}$$

*Proof.* Let $\beta \geq 2$, $\mathbf{u} \in \mathcal{U}_\beta(m) \cap \{0,1\}^\mathbb{N}$ and $\tilde{\mathbf{u}}$ as in Lemma 3.16. Then $\tilde{\mathbf{u}} \in \mathcal{U}_{\tilde{\beta}}(m)$ for all $\tilde{\beta} > \beta$. If $\tilde{\mathbf{u}} \in \mathcal{S}_\infty$, then Lemma 3.10 gives that $\beta \geq 1 + \sqrt{m}$. If $\tilde{\mathbf{u}} = \sigma(\overline{1})$, $\sigma \in S^*$, then $\beta \geq \max\{f_\sigma(m), g_\sigma(m)\}$ by Lemmas 3.11 and 3.12. This implies that $\beta \geq f_\sigma(m)$ if $m \in \left[\mathfrak{m}_{\sigma(0\overline{1})}, \mu_\sigma\right]$, $\beta \geq g_\sigma(m)$ if $m \in \left[\mu_\sigma, \mathfrak{m}_{\sigma(\overline{1})}\right]$, and $\beta \geq 1 + \sqrt{m}$ otherwise. The opposite inequalities are also proved in Lemmas 3.10, 3.11 and 3.12. $\square$

The previous lemmas prove Propositions 3.4 and 3.5 and the main part of Theorem KLP.

*Proof of Theorem 2.* By Lemmas 3.11 and 3.12, $\mathcal{G}(m)$ is differentiable on $(1,2] \setminus \big(\mathfrak{M} \cup \{\mu_\sigma : \sigma \in S^*\}\big)$. By Propositions 3.3 and 3.18, Lemmas 3.11 and 3.12, and the continuity of $\mathcal{G}$ on $(1,2]$, the total variation is

$$\sum_{\sigma \in S^* \setminus S^* \tau_0} \big(\mathcal{G}(\mathfrak{m}_{\sigma(0\overline{1})}) - \mathcal{G}(\mu_\sigma)\big) + \sum_{\sigma \in S^* \setminus S^* \tau_0: \sigma(1) \neq 1} \big(\mathcal{G}(\mathfrak{m}_{\sigma(\overline{1})}) - \mathcal{G}(\mu_\sigma)\big).$$

As $\lim_{m \to 1^+} \mathcal{G}(m) = 2 = \mathcal{G}(2)$, the two sums are equal. For $m \in \big(\mu_\sigma, \mathfrak{m}_{\sigma(\overline{1})}\big)$, $\sigma \in S^*$, $\sigma(1) \neq 1$, we have

$$\frac{1}{\mathcal{G}'(m)} = \frac{m}{\mathcal{G}(m) - 1} - (\mathcal{G}(m) - 1) \sum_{k=1}^\infty \frac{(k+1)u_k^{(\sigma)}}{\mathcal{G}(m)^{k+2}} = 1 - \sum_{k=1}^\infty \left(k - \frac{k+1}{\mathcal{G}(m)}\right) \frac{u_k^{(\sigma)}}{\mathcal{G}(m)^{k+1}}$$

$$> 1 - \sum_{k=1}^\infty \left(k - \frac{k+1}{\mathcal{G}(m)}\right) \frac{1}{\mathcal{G}(m)^{k+1}} = 1 - \frac{1}{\mathcal{G}(m)^2} \geq \frac{3}{4},$$

using that $k - \frac{k+1}{\mathcal{G}(m)} \geq 0$ for all $k \geq 1$. Therefore, we have

$$\sum_{\sigma \in S^* \setminus S^* \tau_0: \sigma(1) \neq 1} \big(\mathcal{G}(\mathfrak{m}_{\sigma(\overline{1})}) - \mathcal{G}(\mu_\sigma)\big) < \frac{4}{3} \sum_{\sigma \in S^* \setminus S^* \tau_0: \sigma(1) \neq 1} (\mathfrak{m}_{\sigma(\overline{1})} - \mu_\sigma) = \frac{4}{3}\left(1 - \sum_{\sigma \in S^* \setminus S^* \tau_0} (\mu_\sigma - \mathfrak{m}_{\sigma(0\overline{1})})\right).$$

We have $\beta_{\tau_h} = 2$ for all $h \geq 0$ since $\tau_h(\overline{1}) = \overline{0^h 1}$ and $\tau_h(0\overline{1}) = 0\overline{0^h 1}$, thus $\mu_{\tau_h} = 2^{h+1}/(2^{h+1} - 1)$, and $\beta = 1 + \sqrt{\mathfrak{m}_{\tau_h(0\overline{1})}}$ satisfies $\beta^{h+3} - 2\beta^{h+2} - 1 = \beta^2 - 2\beta$, hence $\mu_{\tau_0} - \mathfrak{m}_{\tau_0(0\overline{1})} \approx 0.24512$ (and $\mu_\sigma = \mu_{\tau_0}$ when $\sigma$ is the identity, $\mathfrak{m}_{0\overline{1}} = \mathfrak{m}_{\tau_0(0\overline{1})}$), $\mu_{\tau_1} - \mathfrak{m}_{\tau_1(0\overline{1})} \approx 0.05136$. This gives that $\sum_{\sigma \in S^* \setminus S^* \tau_0} (\mu_\sigma - \mathfrak{m}_{\sigma(0\overline{1})}) > 1/4$, thus the total variation is less than 2.

The derivative is unbounded because we have, for all $m \in \big(\mathfrak{m}_{\sigma(0\overline{1})}, \mu_\sigma\big)$ with $\sigma \in \tau_h S^*$,

$$\left|\frac{1}{\mathcal{G}'(m)}\right| = \sum_{k=1}^\infty \frac{k \tilde{u}_k^{(\sigma)}}{\mathcal{G}(m)^{k+1}} \leq \sum_{k=h+1}^\infty \frac{k}{2^{k+1}} = \frac{h+2}{2^{h+1}}. \quad \square$$

*Proof of Theorem 4.* Note that the map $\iota : m \mapsto \frac{m}{m-1}$ is an order-reversing involution on $(1, \infty)$. By Proposition 3.18, we have $m = \iota(\mathcal{G}(m)^h)$ and thus $\mathcal{G}(m) = \iota(m)^{1/h}$ for all $m \in [\mathfrak{m}_{\tau_{h-1}(0\overline{1})}, \mu_{\tau_{h-1}}]$, $h \geq 1$. Moreover, $\iota(m) \geq \mathfrak{m}_{\tau_{h-1}(0\overline{1})}$ is equivalent to $1 + \sqrt{\iota(m)} \geq m^{1/h}$ by Lemma 3.11, and $\mu_{\tau_{h-1}} = 2$ by the proof of Theorem 2. For $2^h \leq m \leq \big(1 + \sqrt{\iota(m)}\big)^h$, we have thus $\mathcal{G}(\iota(m)) = m^{1/h}$. $\square$

### 3.4 Hausdorff dimension of $\mathfrak{M}$

In this section we show that the Hausdorff dimension of $\mathfrak{M}$ is 0.

*Proof of Theorem 3.* It suffices to show that $\dim_H(\mathcal{G}(\mathfrak{M})) = 0$ because $\mathcal{G} : \mathfrak{M} \to \mathbb{R}$ is given by $\mathcal{G}(m) = 1 + \sqrt{m}$, and $1 + \sqrt{m}$ is bi-Lipschitz on the interval $(1, 2]$.

Given $m \in \mathfrak{M}$, we know by Proposition 3.5 that $m = \mathfrak{m}_{\mathbf{u}}$ for some $\mathbf{u} \in \mathcal{S} \setminus \{\overline{1}\}$. Remark 3.2 states that $2\mathbf{u}$ is also the $\beta$-expansion of $\beta = 1 + \sqrt{\mathfrak{m}_{\mathbf{u}}}$. Therefore, for each $n \in \mathbb{N}$ we have

$$1 + \sqrt{\mathfrak{m}_{\mathbf{u}}} \in C_{2u_1 \cdots u_n} := \{\beta > 1 : \text{the } \beta\text{-expansion of } \beta \text{ starts with } 2u_1 \cdots u_n\}.$$

We have $C_{2u_1 \cdots u_n} \subset [2, \infty)$ and, hence, the diameter of $C_{2u_1 \cdots u_n}$ is at most $2^{-n}$, e.g., by a lemma of Schmeling [10, Lemma 4.1].

We now prove $\dim_H(\mathcal{G}(\mathfrak{M})) = 0$ by explicitly constructing a cover. We introduce the set

$$L_n := \{u_1 \cdots u_n \in \{0,1\}^n : u_1 \cdots u_n \text{ is a prefix of an element of } \mathcal{S}\}.$$

For each $n \in \mathbb{N}$ we have

$$\mathcal{G}(\mathfrak{M}) \subset \bigcup_{u_1 \cdots u_n \in L_n} C_{2u_1 \cdots u_n}.$$

So the set $\{C_{2u_1 \cdots u_n} : u_1 \cdots u_n \in L_n\}$ is a cover of $\mathcal{G}(\mathfrak{M})$. Let $s > 0$ be arbitrary and $\mathcal{H}^s(\cdot)$ denote the $s$-dimensional Hausdorff measure. We observe

$$\mathcal{H}^s(\mathcal{G}(\mathfrak{M})) \leq \lim_{n \to \infty} \sum_{u_1 \cdots u_n \in L_n} \text{Diam}(C_{2u_1 \cdots u_n})^s \leq \lim_{n \to \infty} \frac{\#L_n}{2^{ns}}.$$

As was pointed out in Remark 3.1, every element of $\mathcal{S}$ is a Sturmian sequence. Thus it is a consequence of [5, Theorem 2.2.36] that $\#L_n$ grows at most polynomially in $n$. Therefore $\lim_{n \to \infty} \#L_n 2^{-ns} = 0$ and $\dim_H(\mathfrak{M}) \leq s$. Since $s$ is arbitrary we are done. □

## 4 Behaviour at the generalised golden ratio

In this section we discuss the behaviour of the univoque set at the generalised golden ratio. It was observed in [1] that when $\beta = \mathcal{G}(\{0, 1, \ldots, m\})$ for some positive integer $m$, then every $x \in (0, \frac{m}{\beta-1})$ either has a countable infinite of expansions, or a continuum of expansions. In other words $\mathcal{U}_{\mathcal{G}(\{0,1,\ldots,m\})}(\{0, 1, \ldots, m\})$ is still trivial. However, Lemma 3.10 demonstrates that this is not always the case. Indeed the following result is an immediate consequence of this lemma.

**Proposition 4.1.** *If $\mathbf{u} \in \mathcal{S}_\infty$ then $\mathcal{U}_{\mathcal{G}(\mathfrak{m}_{\mathbf{u}})}(\mathfrak{m}_{\mathbf{u}})$ is non-trivial.*

In [9] it was shown that the smallest $\beta \in (1, 2)$ for which an $x$ has precisely two expansions over the alphabet $\{0, 1\}$ was $\beta_2 \approx 1.71064$. As such, there is a small gap between the golden ratio for the alphabet $\{0, 1\}$, and the smallest $\beta$ for which an $x$ has precisely two expansion. As we show below, for certain alphabets it is possible that an $x$ has precisely two expansions at the golden ratio.

**Proposition 4.2.** *For every $\mathbf{u} \in \mathcal{S}_\infty$, the number $\mathfrak{m}_{\mathbf{u}}/\mathcal{G}(\mathfrak{m}_{\mathbf{u}})$ has precisely two expansions in base $\mathcal{G}(\mathfrak{m}_{\mathbf{u}})$ over the alphabet $\{0, 1, \mathfrak{m}_{\mathbf{u}}\}$.*

*Proof.* Let $\beta = \mathcal{G}(\mathfrak{m}_{\mathbf{u}}) = 1 + \sqrt{\mathfrak{m}_{\mathbf{u}}}$ and let $\mathfrak{m}_{\mathbf{u}}/\beta = \sum_{k=1}^{\infty} v_k \beta^{-k}$ be an expansion of $\mathfrak{m}_{\mathbf{u}}/\beta$ over the alphabet $\{0, 1, \mathfrak{m}_{\mathbf{u}}\}$. Since $\mathfrak{m}_{\mathbf{u}} > \frac{\mathfrak{m}_{\mathbf{u}}}{\beta-1}$, we have $v_1 \in \{1, \mathfrak{m}_{\mathbf{u}}\}$, thus $\sum_{k=1}^{\infty} v_{k+1} \beta^{-k}$ equals $\mathfrak{m}_{\mathbf{u}} - 1$ and $0$ respectively. Clearly, $0$ has a unique expansion, and $\mathfrak{m}_{\mathbf{u}} - 1$ has the expansion $u_1 u_2 \cdots$ by (3.1), which is also unique by Lemma 3.10. □

Proposition 4.1 and Proposition 4.2 imply Theorem 5.

**Acknowledgements** The authors are grateful to Vilmos Komornik for posing the questions that lead to this research.

## References

[1] S. Baker, *Generalized golden ratios over integer alphabets,* Integers **14** (2014), Paper No. A15, 28 pp.

[2] Z. Daróczy, I. Katai, *Univoque sequences*, Publ. Math. Debrecen **42** (1993), 397–407.

[3] P. Erdős, I. Joó, V. Komornik, *Characterization of the unique expansions $1 = \sum_{i=1}^{\infty} q^{-n_i}$ and related problems*, Bull. Soc. Math. Fr. **118** (1990), 377–390.

[4] V. Komornik, A.C. Lai, M. Pedicini, *Generalized golden ratios of ternary alphabets,* J. Eur. Math. Soc. (JEMS) **13** (2011), no. 4, 1113–1146.

[5] M. Lothaire, *Algebraic combinatorics on words,* Encyclopedia of Mathematics and its Applications, 90. Cambridge University Press, Cambridge, 2002. xiv+504 pp. ISBN: 0-521-81220-8

[6] W. Parry, *On the $\beta$-expansions of real numbers*, Acta Math. Acad. Sci. Hung. **11** (1960) 401–416.

[7] M. Pedicini, *Greedy expansions and sets with deleted digits,* Theoret. Comput. Sci. **332** (2005), no. 1–3, 313–336.

[8] A. Rényi, *Representations for real numbers and their ergodic properties*, Acta Math. Acad. Sci. Hung. **8** (1957) 477–493.

[9] N. Sidorov, *Expansions in non-integer bases: lower, middle and top orders,* J. Number Th. **129** (2009), 741–754.

[10] J. Schmeling, *Symbolic dynamics for $\beta$-shifts and self-normal numbers,* Ergodic Theory Dynam. Systems **17** (1997), no. 3, 675–694.